\documentclass[11pt]{article}
\usepackage{amssymb,a4,epsfig}
\usepackage{amssymb, latexsym}

\begin{document}

\newcommand{\e}{\varepsilon}
\renewcommand{\a}{\alpha}
\renewcommand{\b}{\beta}
\newcommand{\om}{\omega}
\newcommand{\Om}{\Omega}
\newcommand{\p}{\partial}
\renewcommand{\phi}{\varphi}

\newcommand{\N}{{\mathbb N}}
\newcommand{\R}{{\mathbb R}}
\newcommand{\EX}{{\mathbb E }}
\newcommand{\PX}{{\mathbb P }}

\newcommand{\cF}{{\cal F}}
\newcommand{\cG}{{\cal G}}
\newcommand{\cD}{{\cal D}}
\newcommand{\cO}{{\cal O}}

\newcommand{\ddt}{\frac{d}{dt}}

\newtheorem{theorem}{Theorem}
\newtheorem{lemma}{Lemma}
\newtheorem{remark}{Remark}
\newtheorem{cor}{Corollary}

\title{Predictability of the Burgers dynamics\\ under model uncertainty}

{\tiny
\author{Dirk Bl\"omker$^1$ and Jinqiao Duan$^2$  \\
1. Institut f\"ur Mathematik\\
   RWTH Aachen\\
   52062 Aachen, Germany\\
2.  Department of Applied Mathematics \\
Illinois Institute of Technology \\
Chicago, IL 60616, USA}
}

 \date{June 15, 2006     }

\maketitle

\begin{abstract}

Complex systems may be subject to various uncertainties. A great
effort  has been concentrated on predicting the dynamics under
uncertainty in initial conditions. In the present work, we
consider the well-known Burgers equation with random boundary
forcing or with random body forcing. Our goal is to attempt to
understand the stochastic Burgers dynamics by predicting or
estimating the solution processes in various diagnostic metrics,
such as mean length scale, correlation function and mean energy.
First, for the linearized model, we observe that the important
statistical quantities like mean energy or correlation functions
are the same for the two types of random forcing, even though the
solutions behave very differently. Second, for the full nonlinear
model, we estimate the mean energy for various types of random
body forcing, highlighting the different impact on the overall
dynamics of space-time white noises, trace class white-in-time and
colored-in-space noises, point noises, additive noises or
multiplicative noises.



\end{abstract}

\tableofcontents

\newpage

\section{Introduction}

The   Burgers equation has been used as a simplified prototype
model for hydrodynamics and infinite dimensional systems. It is
often regarded as a one-dimensional Navier-Stokes equation. Our
motivation for considering this  equation  comes from the modeling
of the hydrodynamics and thermodynamics of  the coupled
atmosphere-ocean system. At the  air-sea interface, the atmosphere
and ocean interact through  heat flux  and freshwater flux with a
fair amount of uncertainty \cite{Stocker, DijkBook, Lions92}.
  These translate into random Neumann boundary conditions for
temperature or salinity. The Dirichlet boundary condition is also
appropriate under other physical situations.  The fluctuating wind
stress forcing corresponds to a random body forcing for the fluid
velocity field. The coupled atmosphere-ocean system is quite
complicated and numerical simulation is the usual approach at this
time. In this paper, we consider a simplified
  model for this system, i.e., we consider the
Burgers equation with random Neumann boundary conditions and
random body forcing. Although the stochastic Burgers equation is
widely studied, most work we know are for Dirichlet boundary
conditions or periodic boundary conditions \cite{DaPratoB, E,
Chueshov, Chueshov2}. The reference \cite{Titi} studied the
control of deterministic Burgers equation with Neumann boundary
conditions.

We consider the stochastic Burgers equation with
boundary forcing on the interval $[0,L]$
\begin{equation}             \label{burg-bound}
  \partial_t u +u\cdot \partial_x u =\nu \partial_x^2 u
\end{equation}
\begin{equation}\label{burg-bound2}
  \partial_x u(\cdot,0) = \alpha \eta
  \qquad \partial_x u(\cdot,L) =0.
\end{equation}
Here $\alpha>0$ denotes the noise strength and
$\eta$ is white noise, i.e., $\eta$ is a generalized
Gaussian process with $\EX\eta(t)=0$ and
$\EX \eta(t)\eta(s)=\delta(t-s)$.
The restriction to noise on the left
boundary is only for simplicity. Analogous results
will be true, if forces act on both sides of the domain.

We will see that boundary forcing coincides with point forcing at the
boundary. Thus we also look at point forcing. As a simple example for point forcing,
we consider
\begin{equation}                \label{burg-point}
  \partial_t v +v\cdot \partial_x v =\nu \partial_x^2 v +\alpha\delta_0 \eta
\end{equation}
\begin{equation}
 u(\cdot,-L) =  u(\cdot,L) =0\;,
\end{equation}
where $\delta_0$ is the Delta-distribution.

We will compare solutions of (\ref{burg-bound}) and
(\ref{burg-point}) with solutions of
 the stochastic Burgers equation with
body forcing.
\begin{equation}                \label{burg-body}
  \partial_t v +v\cdot \partial_x v =\nu \partial_x^2 v +\sigma\xi
\end{equation}
either subject to Dirichlet or  Neumann boundary conditions.
Here the noise strength is denoted by $\sigma>0$ and $\xi$ is
space-time white noise. I.e., $\xi$ is a generalized Gaussian
process with $\EX\xi(t,x)=0$ and $\EX
\xi(t,x)\xi(s,y)=\delta(t-s)\delta(x-y)$. We will also consider
trace class body noise, i.e., noise that is white in time but
colored in space.

For the linearized equations, we will compare
statistical quantities of both solutions, which are frequently used.
One of them is the
{\em mean energy}
\begin{equation}                \label{m-e}
  \frac1L \int_0^L  \EX[u(t,x)-\overline{u}(t)]^2 dx,
\end{equation}
where
$$\overline{u}(t)= \frac1L \int_0^L u(t,x) dx.
$$
Another important quantity, which gives information about
the characteristic size of pattern, is the
{\em mean correlation function}
\begin{equation}                \label{m-cf}
  C(t,r):=
  \frac1L \int_0^L \EX [u(t,x)-\overline{u}(t)]\cdot
  [u(t,x+r)-\overline{u}(t)] dx,
\end{equation}
which is usually averaged over all points $r$ with a given distance from $0$.
We obtain the   {\em averaged mean correlation function}
\begin{equation}                \label{m-acf}
\hat{C}(t,r)=\frac12[ C(t,r)+ C(t,-r)]
\end{equation}
where we employ the canonical odd and 2L-periodic
extension of  $u$ in order to define
$C(t,r)$ for any $r\in \R$.

For the  linearized equation the {\em main result} states that
mean energy and averaged mean correlation function are the same
for solutions of (\ref{burg-bound}) and (\ref{burg-body}).
Nevertheless the solutions behave completely different.
Furthermore, we give some qualitative properties like, for
instance, the typical pattern size. This should carry over to  a
transient regime (i.e., small times) for the corresponding
nonlinear equations.

For the full nonlinear Burgers model, we estimate the mean energy
for various types of random body forcing, highlighting the
different impact on the overall dynamics of space-time white
noises, trace class white-in-time and colored-in-space noises,
point noises, additive noises or multiplicative noises.

In the following, we discuss linear dynamics in \S 2 and nonlinear
dynamics in \S 3.


\section{Linear Theory}

Define $$ A=\nu\partial_x^2$$ with $$ D(A)=\{w\in H^2([0,L]) \ : \
\partial_xw(0)=0,\  \partial_xw(L)=0 \}$$ It is well-known (cf.
e.g. \cite{Courant}) that $A$ has an orthonormal base of
eigenfunctions $\{e_k\}_{k\in\N_0}$ in $L^2([0,L])$ with
corresponding eigenvalues $\{\lambda_k\}_{k\in\N_0}$. In our
situation $e_0(x)=1/\sqrt{L}$, $e_k(x)=\sqrt{2/L}\cdot\cos(\pi k
x/L)$ for $k\in\N$, and $\lambda_k=-\nu (k\pi/L)^2$. Moreover $A$
generates an analytic semigroup $\{e^{tA}\}_{t\ge0}.$ (cf. e.g.
\cite{Luna}).

In fact, $e^{tA}v_0$ is the solution of the following
evolution problem
\begin{equation}
  \partial_t v  = A v, \qquad
  \partial_x v(\cdot,0) = \partial_x v(\cdot,L) =0, \qquad v(0, x)=v_0(x).
\end{equation}
The solution is
\begin{equation}    \label{semigroup}
 e^{tA}v_0(x) := v(t,x) = \sum_{k=0}^{\infty} <v_0, e_k>e^{\lambda_k t} e_k,\quad
 t>0, \; 0<x<L,
\end{equation}
where $<\cdot, \cdot>$ is the usual scalar product in $L^2(0, L)$.

We now consider
the following linearized problems. First
\begin{equation}             \label{bound-lin}
  \partial_t u  =A u, \qquad
  \partial_x u(\cdot,0) = \alpha \partial_t\beta,
  \quad \partial_x u(\cdot,L) =0.
\end{equation}
Here   the white noise  $\eta$ is given by the  generalized
derivative of a standard Brownian motion (cf. e.g.
\cite{Arnold2}), and $\alpha$ is the noise intensity.

Secondly,
\begin{equation}                \label{body-lin}
  \partial_t v  = A v + \sigma \partial_t W,
\qquad
  \partial_x v(\cdot,0) = 0,
  \quad \partial_x v(\cdot,L) =0,
\end{equation}
where the space-time white noise is given by the generalized
derivative of an $Id$-Wiener process. Namely, $W(t)=
\sum_{k\in\N_0}\beta_k(t)e_k$, where $\{ \beta_k \}_{k\in\N}$ is a
family of independent standard Brownian motions, and $\sigma$ is
the noise intensity.

It is   known (cf. e.g. \cite{DaPrato}) that (\ref{body-lin}) has
a unique weak solution given by the stochastic convolution (taking
initial condition to be zero)
\begin{equation}
 W_A(t)
   =\sigma\cdot \int_0^t e^{(t-\tau)A}dW(\tau)
   =\sigma\cdot \sum_{k\in\N_0}
          \int_0^t e^{(t-\tau)\lambda_k} d\beta_k(\tau)e_k.
\end{equation}
We define the Neumann map $\cD$ by
$$(1-A)\cD\gamma=0,
    \quad  \partial_x\cD\gamma(0)=\gamma,\
     \partial_x\cD\gamma(L)=0
$$
for any $\gamma\in\R$. It is  known   that $\cD:\R\mapsto
H^2([0,L])$ is a continuous linear operator. In fact, we have
explicit expression for this linear operator

\begin{eqnarray} \label{Neumann}
 \cD (\gamma) = \frac{e^x+e^{2L}e^{-x}}{1-e^{2L}} \gamma.
\end{eqnarray}
        From \cite{DaPrato2} or \cite{DaPratoA} we immediately obtain,
that (\ref{bound-lin}) has a unique weak solution (taking initial
condition to be zero)
\begin{equation}
  Z(t)
= (1-A)\int_0^t e^{(t-\tau)A}\cD \alpha d\beta(\tau).
\end{equation}
In the next section we derive explicit formulas for $Z$ in term of Fourier series.

%
\subsection{Mean Energy}

To obtain the Fourier series expansion for $Z,$ consider
for $e\in D(A)$ and $\gamma\in\R$
\begin{eqnarray}
<\cD(\gamma), (1-A)e>_{L^2([0,L])}
&=& <\cD(\gamma),e>-\int_0^L \cD(\gamma)\cdot e_{xx} dx
\nonumber\\
&=& <\cD(\gamma),e>-\int_0^L \cD(\gamma)_{xx} \cdot e dx
     +\left. \cD(\gamma)_{x} \cdot e\right|_{x=0}^{x=L}
\nonumber\\
&=& - \gamma e(0),
\end{eqnarray}
by the definition of $\cD$. Hence,
\begin{eqnarray}
<Z(t), e_k>
&=& <\int_0^t e^{(t-\tau)A}\cD \alpha d\beta(\tau),
      (1-A)e_k \nonumber>\\
&=& \int_0^t e^{(t-\tau)\lambda_k}< \cD \alpha d\beta(\tau),
      (1-A)e_k >\nonumber\\
&=&\alpha e_k(0)\cdot\int_0^t e^{(t-\tau)\lambda_k} d\beta(\tau).
\end{eqnarray}
We now obtain
\begin{equation}
  Z(t)=\alpha\cdot\sum_{k\in\N_0}
  e_k(0)\int_0^t e^{(t-\tau)\lambda_k} d\beta(\tau)e_k.
\end{equation}
Finally,
\begin{equation}
 Z(t)=\alpha e_1(0) \cdot \sum_{k\in\N}
  \int_0^t e^{(t-\tau)\lambda_k} d\beta(\tau)e_k
  +\alpha e^2_0(0)\beta(t)
\end{equation}
and
\begin{equation}
W_A(t)=\sigma\cdot \sum_{k\in\N}
  \int_0^t e^{(t-\tau)\lambda_k} d\beta_k(\tau)e_k
  +\sigma \beta_0(t)\;.
\end{equation}
If we now choose $\sigma=\alpha e_1(0)$, we readily obtain that
$$
\EX\|Z(t)-\overline{Z}(t)\|^2
=\sigma^2 \cdot \sum_{k\in\N}\int_0^t e^{2\tau\lambda_k} d\tau
=\EX\|W_A(t)-\overline{W_A}(t)\|^2\;,
$$
where $\| \cdot \|$ is the norm in $L^2([0,T]).$
Hence, the mean energy in both cases is given by
$ \sigma^2L^{-1} \sum_{k\in\N}\int_0^t e^{2\tau\lambda_k} d\tau$.

For the mean energy we can prove the following theorem,
which is similar to the results of \cite{DB-rough} and \cite{DB-rough2}.
\begin{theorem}
\label{thm:mebo}
Fix $\sigma^2=\alpha^2/L$, then
the mean energy $C_Z(t,0)=C_{W_A}(t,0)$ behaves like
$C_1(\alpha^2/L)\sqrt{t/\nu}$ for $t\ll L^2/\nu$,
and like $C_2\alpha^2/\nu$ for $t\gg L^2/\nu$.
\end{theorem}
The main difference to body forcing is the scaling in the length-scale $L$.
The long-time scaling is independent of $L$, while the transient scaling is.
 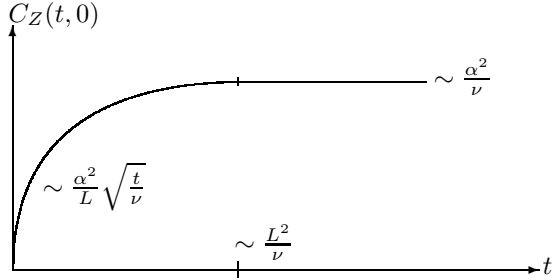
\begin{figure}[h]
 \begin{center}
\setlength{\unitlength}{0.5mm}
\begin{picture}(140,70) 
      (0,0) 
   \put(0,0){\vector(1,0){140}}
   \put(0,0){\vector(0,1){65}}
   \qbezier(0,0)(0,49)(60,50)
   \put(60,50){\line(1,0){50}}
   \small
   \put(60,2){\line(0,-1){4}}
   \put(59,5){$\sim \frac{L^2}{\nu}$}
   \put(141,-1){$t$}
   \put(-1,66){$C_Z(t,0)$}
   \put(112,49){$\sim  \frac{\alpha^2}{\nu}$}
   \put(8,20){$\sim  \frac{\alpha^2}{L}\sqrt{\frac{t}{\nu}}$}
   \put(60,51){\line(0,-1){2}}
\end{picture}
 \end{center}
 \vspace*{-5mm}
 \caption{The scaling of the mean energy for boundary forcing.}
 \end{figure}

%
\subsection{Correlation Function}
\label{sec:corfu}

To obtain results for the correlation function, we think of $Z$ and $W_A$ to be
periodic on $[-L,L]$, and symmetric w.r.t. $0$.
I.e., we choose the standard $2L$-periodic extension respecting the Neumann boundary conditions on $[0,L]$.
To be more precise, we extend $Z$ and $W_A$ in a Fourier series in the basis $e_k$, which we then consider
to be defined on  whole $\R$.

We consider firstly for $k,l\not=0$
$$
  \int_0^L e_k(x)e_l(x+r) dx
    =\left\{
\renewcommand{\arraystretch}{1.5}
\begin{array}{ccc}
       L e_k(0)e_l(0)\frac{l((-1)^{k+l}-1)}{\pi(l^2-k^2)} \sin(\pi l r/L) &:& k \not= l\\[1mm]
       L e_k(0)^2\cos(\pi k r/L) &:& k=l
     \end{array}\right..
$$
Now relying on the independence of the Brownian motions, it is straightforward to verify
\begin{eqnarray}    \label{repr-CWA}
C_{W_A}(t,r)
& =& \frac1L \EX<W_A(t,x)-\overline{W_A}(t),W_A(t,x+r)-\overline{W_A}(t)>
     \nonumber\\
& =& \frac{\alpha^2e_1^2(0)}{L}\cdot
      \sum_{k\in\N}\int_0^t e^{2\tau\lambda_k}d\tau
         \cos(\pi k r/L)\;,
\end{eqnarray}
as $\alpha^2e_1^2(0)=\sigma^2$.
Furthermore,
\begin{eqnarray}    \label{repr-CZ}
 C_{Z}(t,r)
& =& \frac1L \EX<Z(t,x)-\overline{Z}(t),Z(t,x+r)-\overline{Z}(t)>\\
& =& C_{W_A}(t,r)+ \frac{e_1^2(0)}{\pi} \sum_{k,l=1 \atop k\not=l}^\infty \int_0^t e^{\tau(\lambda_k+\lambda_l)} d\tau  \frac{l ((-1)^{k+l}-1)}{(l^2-k^2)} \sin(\pi l r/L)\;.\nonumber
\end{eqnarray}
Obviously, $C_Z$ and $C_{W_A}$ do not coincide, but
let us now look at the averaged correlation function
$$\hat{C}(t,r)=\frac12[C(t,r)+C(t,-r)]\;.
$$
Then it is obvious that
\begin{equation}
\label{e:avCZ}
\hat{C}_{W_A}(t,r)=C_{W_A}(t,r)=\hat{C}_{Z}(t,r)\not=C_Z(t,r)\;.
\end{equation}
Now
\begin{theorem} \label{lin-co}
For $\alpha^2e_1^2(0)=\sigma^2$ the mean energy
and the averaged mean correlation functions $\hat{C}_{W_A}$ and
$\hat{C}_{Z}$ for $Z$ and $W_A$
coincide for any $t\ge0$.
\end{theorem}

\begin{figure}[h]
\begin{center}
\includegraphics[angle=0,width=9cm]{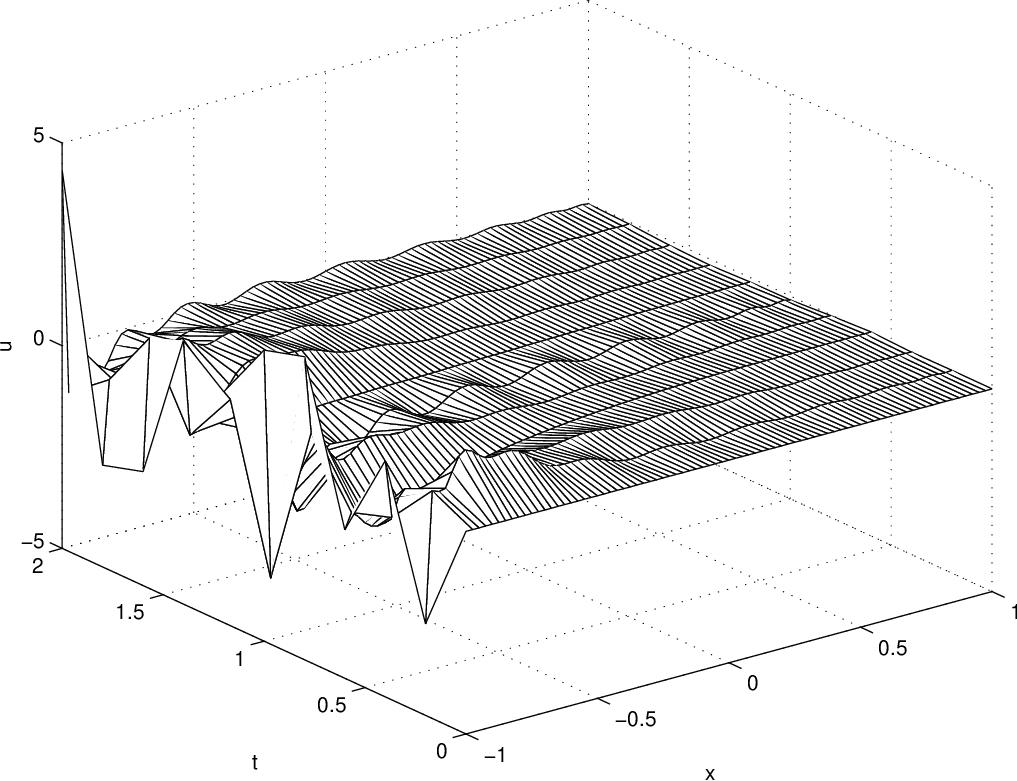}
\end{center}
\caption{\emph{Random boundary condition}: One realization of the
solution of the equation (\ref{bound-lin}) for $L=1$, $\nu=1$,
$\alpha=1$  and initial condition $u(x,0)=0$.}\label{RandomBC}
\end{figure}

\begin{figure}[h]
\begin{center}
\includegraphics[angle=0,width=9cm]{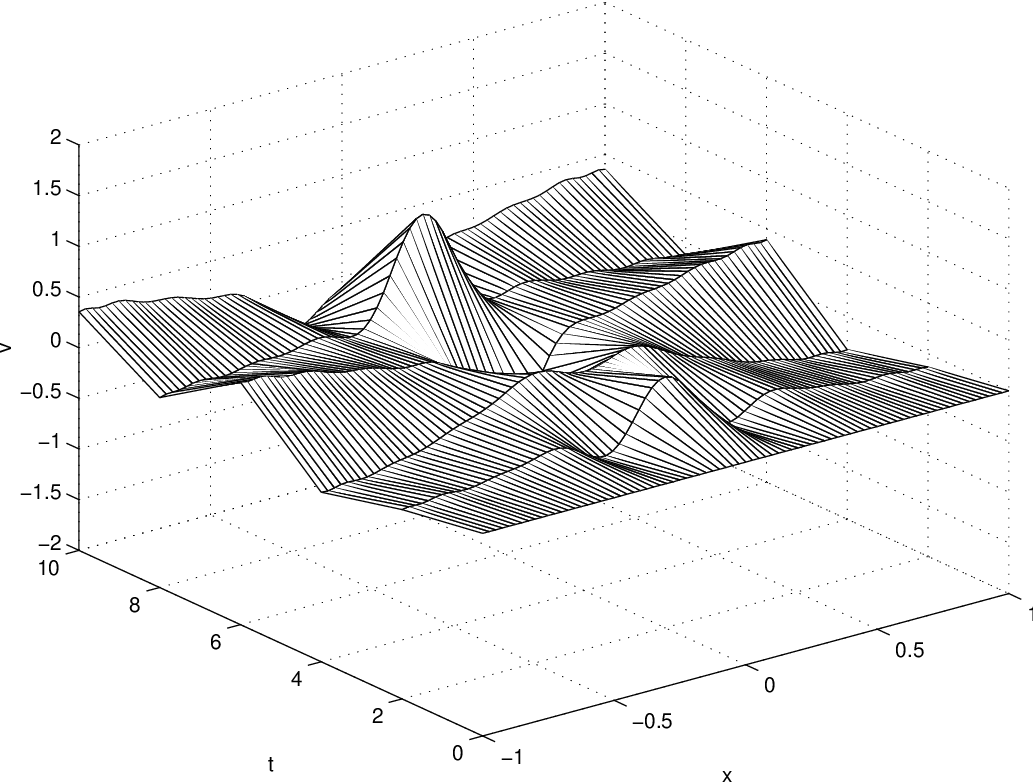}
\end{center}
\caption{\emph{Random body forcing}: One realization of the
solution of the equation (\ref{body-lin}) for $L=1$,  $\nu=1$,
$\sigma=1$ and initial condition $v(x,0)=0$.}\label{BodyForcing}
\end{figure}

This is somewhat surprising, as realizations of $Z$ and $W_A$
behave completely different, when the condition
$\alpha^2e_1^2(0)=\sigma^2$ is satisfied. See e.g. Figure
\ref{RandomBC} and Figure \ref{BodyForcing}.

It is even more surprising, as the scaling behavior of quantities
like mean energy and mean correlation functions are an important tool in
applied science, which for example is used to
determine the size of characteristic length scales and the
universality class the model belongs to. Here both linear models
lie in the same class, although their behavior differs completely.

The scaling behavior with respect to $L$ and $t$ of the mean energy
can be described using the results of \cite{DB-rough},
where the mean surface width for very general models was discussed.
Therefore we focus on the scaling properties of the mean
correlation function. Here we also want to investigate the
dependence on $\alpha$ and $\nu$.

First we consider the scaling properties of the correlation function
$\hat{C}_Z(t,r)$ or $C_{W_A}$, as given in (\ref{repr-CZ}).
We are especially interested in the smallest zero of the function,
which gives information about characteristic length scales or pattern sizes.
For this, we use the {\em normalized correlation function}.
\begin{equation}\label{def-rhoZ}
  \rho_Z(t,r)=\frac{\hat{C}_Z(t,r)}{\hat{C}_Z(t,0)}.
\end{equation}
Note that $C_Z(t,0)$ is the mean energy
and the maximum of $r\mapsto \hat{C}_Z(t,r)$.

We begin with some technical results.
For any continuously differentiable and integrable function
$f:\R^+\to\R$ we obtain using the mean value theorem
\begin{equation}
\left|\int_0^\infty f(x)dx-\sum_{k=1}^\infty f(k) \right|
\le \sum_{k=1}^\infty\sup_{\eta\in(k-1,k)}|f'(\eta)|.
\end{equation}
For $f(k):=e^{-2\tau\nu k^2 \pi^2/L^2}\cos(k\pi r/L)$ it is easy to verify that
$$|f'(k)|\le  e^{-\tau\nu k^2 \pi^2/L^2}\cdot
\frac\pi{L}\cdot
\left[r +\frac{4\sqrt{\tau\nu}}{\sqrt{2 e}}\right],
$$
where we used that $x^se^{-x^2 \alpha} \le (2\alpha e)^{-1/2}$ for
any $x,\alpha \ge0$. Hence,
\begin{eqnarray}
\sum_{k=1}^\infty\sup_{\eta\in(k-1,k)}|f'(\eta)|
&\le&
\sum_{k=1}^\infty  e^{-\tau\nu (k-1)^2 \pi^2/L^2}\cdot
\frac\pi{L}\cdot[r +\frac{4}{\sqrt{2 e}}\sqrt{\tau\nu}]
\nonumber\\ &\le&
\frac\pi{L}[r +\frac{4}{\sqrt{2 e}}\sqrt{\tau\nu}]
(1 + \int_0^\infty e^{-\tau\nu k^2 \pi^2/L^2} dk)
\nonumber\\ &=&
\frac\pi{L}[r +\frac{4}{\sqrt{2 e}}\sqrt{\tau\nu}]
(1 + \frac{L}{\pi\sqrt{\tau\nu}}\cdot \int_0^\infty e^{- k^2} dk)
\end{eqnarray}
and
\begin{equation}
  \int_0^t \sum_{k=1}^\infty\sup_{\eta\in(k-1,k)}|f'(\eta)| d\tau
\le C\frac{1}{L} \sqrt{\frac{t}{\nu}}[r+\sqrt{t\nu}][L+\sqrt{t\nu}].
\end{equation}
Moreover,
\begin{eqnarray}\nonumber
\frac1L\int_0^t\int_0^\infty f(k)dkd\tau
&=& \frac1\pi \int_0^\infty \frac{1-e^{-2t\nu k^2}}{2\nu k^2}
\cos(kr)dk\\
&=& \frac1{\pi}\sqrt{\frac{t}\nu}\cdot G(\frac{r}{\sqrt{\nu t}}),
\end{eqnarray}
with $G(x):=\int_0^\infty  \frac{1-e^{-2 k^2}}{2 k^2} \cos(k x)dk$.

 \begin{figure}[h]
 \begin{center}
 \epsfxsize 0.7\hsize \epsfbox{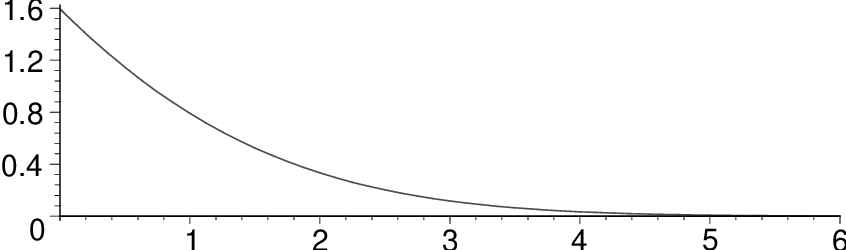}
 \end{center}
 \vspace*{-5mm}
 \caption{A sketch of G}
 \end{figure}

Using  (\ref{e:avCZ}) we immediately obtain
$$ C_{Z}(t,r)= \frac{\alpha^2}{L}
\frac1{2\pi} \sqrt{\frac{t}\nu}\cdot G(\frac{r}{\sqrt{\nu t}})
+\cO\left( \frac{\alpha^2\sqrt{t}}{L^3\sqrt{\nu}}[r+\sqrt{t\nu}][L+\sqrt{t\nu}] \right)\;.
$$
Note that the approximation with $G$ is not $L$-periodic in $r$,
while $C_{Z}(t,r)$ is. The solution is that the error term is $\cO(1)$ for $r$ near $L$.

For the normalized correlation function we deduce
$$ \rho_Z(t,r) := \frac{\hat{C}_Z(t,r)}{\hat{C}_Z(t,0)}
= \frac{G(\frac{r}{\sqrt{\nu t}})
+\cO\left(\frac1{L^2}[r+\sqrt{t\nu}][L+\sqrt{t\nu}]\right)}
{G(0)+\cO\left(\frac1{L^2}[\sqrt{t\nu}][L+\sqrt{t\nu}] \right)}.
$$
           From the properties of $G$ we infer the following:
\begin{theorem}
Given $\delta\in(0,1)$ and sufficiently small $\epsilon_2>0,$
there exists some $\epsilon_1>0$ and three constants $0<C_1<C_2<C_3$
depending only on $\delta$ and $\epsilon_2$ such that
for $t <\epsilon_1 L^2/\nu$ the following holds:
$$\rho_Z(t,r)\ge \delta
\quad\mathrm{\it for}\quad r\in[0,C_1\sqrt{t\nu}]
$$
and
$$|\rho_Z(t,r)|<\epsilon_2
\quad\mathrm{\it for}\quad r\in[C_2\sqrt{t\nu},C_3\sqrt{t\nu}].$$
\end{theorem}
Note that we did not show that the correlation function has a zero,
but it is arbitrary small in a point $r_\epsilon\approx \sqrt{t\nu}$.
Therefor the theorem says that the typical length-scale is $\sqrt{t\nu}$,
at least for times $t\ll L^2/\nu$.

For $t\to\infty$ we immediately obtain that
$$ \hat{C}_Z(\infty,r)
= \frac{\alpha^2}{\pi^2\nu}
\sum_{k=1}^\infty \frac1{k^2} \cos(k\pi r/L)
=: \frac{\alpha^2}{\nu} F(r/L)
$$
and
$$|\hat{C}_Z(t,r)-\hat{C}_Z(\infty,r)|
\le \frac{\alpha^2}{\pi^2\nu} e^{-2t\nu\pi^2/L^2}\sum_{k=1}^\infty\frac1{k^2}.
$$

\begin{figure}[h]
\begin{center}
\epsfxsize 0.6\hsize \epsfysize 0.4\hsize \epsfbox{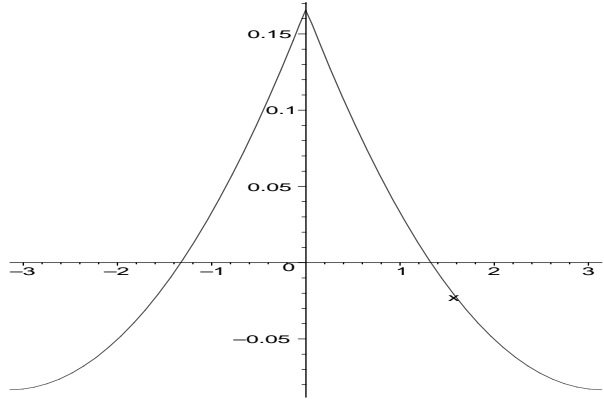}
\end{center}
\vspace*{-5mm}
\caption{A sketch of F}
\end{figure}

We can look for the explicit representation of $F$, which is a
$2$-periodic function, and compute explicitly the zero, but all we need from $F$ is, that for a
given small enough $\delta>0$ there is a $x_\delta>0$ such that
$F>\delta$ on $[0,x_\delta]$. Moreover, there is some $x_0$ such
that $F(x_0)=0$.

Consider the normalized correlation function
$$ \rho_Z(t,r) = \frac{\hat{C}_Z(t,r)}{\hat{C}_Z(t,0)}
=\frac{F(r/L)+\cO(e^{-2t\nu\pi^2/L^2})}{F(0)+\cO(e^{-2t\nu\pi^2/L^2})}.
$$
Assume that $t\nu\gg L^2$ (i.e., there is some small  $\epsilon>0$
such that $\epsilon t\nu> L^2$). Now,
$$  \rho_Z(t,x_0 L)=\cO(e^{-2t\nu\pi^2/L^2})$$
and
$$ \rho_Z(t,x L)
\ge \frac\delta{F(0)} + \cO(e^{-2t\nu\pi^2/L^2})>0$$
for any $x<x_\delta$.

So for $t\gg L^2/\nu$ the first zero of $\rho_Z$ should be
of order $L$. A more precise formulation is:

\begin{theorem}
Given $\delta\in (0,0.8)$  and $\delta\gg\epsilon_2>0$,
there exists some $\epsilon_1>0$, a constant $C>0$,
and a point $x_o>0$
depending only on $\delta$ and $\epsilon_2$ such that
for $t > L^2/ (\nu\epsilon_1)$ we obtain the following:
$$\rho_Z(t,r)\ge \delta
\quad\mathrm{\it for}\quad r\in[0,C L]
$$
and
$$|\rho_Z(t,x_0 L)|<\epsilon_2.$$
\end{theorem}
Thus the theorem tells us that for $t\gg L^2/\nu$, the typical length scale is of order $L$,
which is the size of system. This result is true for both boundary and body forcing.


\section{Nonlinear theory}

For the nonlinear results we leave the  setting of boundary forcing.
Mainly, due to the lack of a-priori estimates. Usually, for Neumann boundary conditions
one relies on the maximum principle to bound solutions,
but the solution for boundary forcing is quite rough. Therefore, we hardly
get useful results. Only, transient bounds for small times are possible to establish.
For the next sections, we focus first on body forcing and later on point forcing.
We will see later that boundary forcing is actually just a point forcing in a point at the boundary.

The main results of this sections are uniform bounds on the energy and thus on the
correlation function $C$, as
$|C(t,r)|\le C(t,0)$, and $C(t,0)$ is the energy.
Furthermore, we show that for $t\to 0$ the linear regime dominates.
In \cite{Du-Bl-Wa} also H"older-continuity for the mean energy was shown for
a quasigeostrophic model.
We conclude this section by a qualitative discussion on upper bounds for the energy
using additive and multiplicative trace-class noise.

\subsection{Body forcing - Mean energy bounds}

Here we provide bounds on the mean energy for the body forcing case.
We consider additive space-time white noise case first,
and show that the mean energy and thus
 the correlation function is uniformly bounded in time.
This result is known (cf. \cite{Le-Nu-Pe:00}) for Burgers equation
using the celebrated Cole-Hopf transformation, but we provide here
a simple  proof for completeness. Furthermore, our proof is based
on energy estimates and it is easily adapted to other types of
equations and additional terms in the equation. In contrast to
that Cole-Hopf transformation is strictly limited to the standard
Burgers equation.

For a long time for space-time white noise only uniform bounds for logarithmic moments were known.
See \cite[Lemma 14.4.1]{DaPrato2} or \cite{DP-Ga:95}.
In \cite{Le-Nu-Pe:00} the transformation to a stochastic heat equation via the celebrated Cole-Hopf
transformation was used to study finiteness of moments. Here we rely on a much simpler tool, which can also
be applied to other equations. See for instance \cite{Du-Bl-Wa} for a quasigeostrophic model,
where our analysis would apply, too.

Consider
\begin{equation}
\label{e:burgbodfo}
  \partial_t u+ u\cdot \partial_x u  =\nu \partial_x^2 u +\sigma \partial_t W\;,
\end{equation}
\begin{equation}
 u(\cdot,-L) = u(\cdot,L) =0, \qquad u(x,0)=u_0(x).
\end{equation}
Here $W$ is a $Q$-Wiener process with a continuous operator $Q\in {\cal L}(L^2)$.
Thus $W$ might be cylindrical, and we include the case of space-time white noise.

Using the semigroup $e^tA$ the
solution for this system is (see \cite{DaPrato2} or
\cite{DaPratoB}):
\begin{equation} \label{mildsolu}
u(t) =  e^{tA}u_0 -\int_0^t e^{(t-\tau)A}(\lambda \Phi_\lambda(\tau)+\frac12\partial_x u(\tau,x)^2) d\tau +  \Phi_\lambda(t)\;,
\end{equation}
where for some $\lambda\ge0$ fixed later
$$ \Phi_\lambda(t)= \sigma \alpha\int_0^t e^{(t-\tau)(A-\lambda)} dW(\tau)
$$
solves
$$ \partial_t \Phi = \nu \partial_x^2 \Phi -\lambda \Phi+ \sigma \partial_t W
$$
subject to
$$  \Phi(\cdot,-L) = \Phi(\cdot,L) =0, \qquad \Phi(x,0)=0\;.
$$

Our main result is now:
\begin{theorem}
\label{thm:e:mebobofo}
Consider  initial conditions  $u_0$ with $ E\|u_0\|^2 <\infty$,
which are independent of the Wiener process $W$ (e.g. deterministic).
Then the mean energy of the solution of (\ref{mildsolu})
is uniformly bounded in time. I.e.,
$$ \sup_{t\ge 0}\EX \|u(t)-\overline{u}(t)\|^2 <\infty\;.
$$
\end{theorem}
\begin{remark}
Actually, we prove that $ \sup_{t\ge 0}\EX \|u(t)\|^2 <\infty$.
The main problem in the proof is that after applying
Gronwall-type estimates we end up with terms $\EX \exp\{\int_0^t\|\Phi_\lambda(s)\|^2_{L^\infty}\}.$
This might blow up in finite time, as second order exponential moments of the Gaussian
may fail to exist, if $t$ is too large. This is why we introduced artificially
additional dissipation in the equation for $\Phi_\lambda$,
in order to get exponential moments small.
\end{remark}
For the proof of Theorem \ref{thm:e:mebobofo} define
\begin{equation}
 \label{e:defv}
v(t)=u(t)-\Phi_\lambda(t)\quad \mathrm{for}\ t\ge0,\ \lambda\ge0\;.
\end{equation}
We see that $v$ is a weak
solution of
\begin{equation}
\label{e:veq}
  \partial_t v+ \frac12\partial_x (v + \Phi_\lambda)^2  = \nu \partial_x^2 v+\lambda \Phi_\lambda
\end{equation}
\begin{equation}
 v(\cdot,-L) = v(\cdot,L) =0, \qquad v(x,0)=u_0(x).
\end{equation}
The following calculation is now only formal, but it can easily made
rigorous using for instance spectral Galerkin approximations. Taking the scalar product in (\ref{e:veq})
yields
\begin{eqnarray*}
\frac12\partial_t \|v\|^2
&=&   -\|v_x\|^2 + \int_{-L}^L (v+\Phi_\lambda)^2 v_x \;dx +  \int_{-L}^L \Phi_\lambda v\; dx\\
&\le&
 -\|v_x\|^2+\|\Phi_\lambda\|_{H^{-1}}\|v_x\|
    +|v_x\|\|\Phi_\lambda\|^2_{L^4}+ 2 \|v_x\|\|\Phi_\lambda\|_{L^\infty}\|v\|\\
&\le& -\frac12c_p^2 \|v\|^2 + 4\|v\|^2\|\Phi_\lambda\|_\infty^2
    + 2\lambda^2\|\Phi_\lambda\|_{L^4}^4+2\|\Phi_\lambda\|_{H^{-1}}^2\;,
\end{eqnarray*}
where we used Young inequality ($ab\le\frac18a^2+2b^2$), and Poincare-inequality $\|v\|\le c_p\|v_x\|$.
Now, from Gronwall-type inequalities
\begin{eqnarray}\label{e:apri}
\|v(t)\|^2
&\le& e^{ -c_p^2t +8\int_0^t\|\Phi_\lambda\|^2_\infty d\tau } \|u(0)\|^2\\
&& + \int_0^t e^{ -c_p^2(t-s) +8\int_s^t\|\Phi_\lambda\|_\infty^2 d\tau }
4(\lambda^2\|\Phi_\lambda\|_{L^4}^4+ \|\Phi_\lambda\|_{H^{-1}}^2)ds \nonumber
\end{eqnarray}
Now we use the following lemma,
which is easily proved by Fernique's theorem,
if we consider $\Phi_\lambda$ as a Gaussian in $L^2([0,t_0],L^\infty)$.
\begin{lemma}
\label{lem:Fernique}
Fix $K>0$ and $t_0>0$, then there is a $\lambda_0$ such that
$$ \sup_{t\in[0,t_0]}\EX \exp\{ 16\int_0^t \|\Phi_\lambda(s)\|^2_{L^\infty}ds \} \le  K^2 $$
for all $\lambda \ge \lambda_0$.
\end{lemma}

Furthermore, we use that all moments of $\|\Phi_\lambda\|_{L^\infty}$
and $\|\Phi_\lambda\|_{H^{-1}}$ are uniformly bounded in time.
This is easily proven, using for instance the celebrated factorization method.

Now we first fix $K>0$, and then $t_0$ such that $e^{ -c_p^2t}K <\frac14$.
This yields for $t\in[0,t_0]$ and  $\lambda$ sufficiently large
$$
\EX\|v(t)\|^2 \le e^{ -c_p^2t} K \EX\|u(0)\|^2 + 4K\int_0^t e^{ -c_p^2(t-s)}
\Big(\EX(\lambda^2\|\Phi_\lambda\|_{L^4}^4+ \|\Phi_\lambda\|_{H^{-1}}^2)^2\Big)^{1/2}ds\;,
$$
using H\"older, Lemma \ref{lem:Fernique}, and the independence of $u(0)$ from $\Phi_\lambda$.
We now find a constant $C$ depending on $t_0$ and $K$ such that
$$
\sup_{t\in[0,t_0]}\EX\|v(t)\|^2 \le  K \EX \|u(0)\|^2 + C
\qquad \mathrm{and} \qquad
\EX\|v(t_0)\|^2 \le  \frac14 \EX \|u(0)\|^2 + C
$$
Using $\EX\|u(t)\|^2 \le 2\EX\|v(t)\|^2+2\EX\|\Phi_\lambda(t)\|^2$
yields for a different constant $C$
$$
\sup_{t\in[0,t_0]}\EX\|u(t)\|^2 \le  K \EX \|u(0)\|^2 + C
\qquad \mathrm{and} \qquad
\EX\|u(t_0)\|^2 \le  \frac12 \EX \|u(0)\|^2 + C
$$
Now we repeat the argument for $k\in\mathbb{N}$ by defining $v(t)=u(kt_0+t)-\tilde\Phi_\lambda(t)$,
where $\tilde\Phi_\lambda(t)$ has the same distribution than $\Phi_\lambda(t)$ due to a time shift of
the Brownian motion. Now $v$ solves again (\ref{e:veq}) with initial condition $u(kt_0)$.
Note that by construction $u(kt_0)$ is independent of $\tilde\Phi_\lambda$.

Repeating the arguments as before yields  for $k\in\mathbb{N}_0$
$$
\sup_{t\in[0,t_0]}\EX\|u(t+kt_0)\|^2 \le  K \EX \|u(kt_0)\|^2 + C
$$
and
$$
\EX\|u((k+1)t_0)\|^2 \le  \frac12 \EX \|u(kt_0)\|^2 + C\;.
$$
Now the following lemma, which is a trivial statement
on discrete dynamical systems, finishes the proof.
\begin{lemma}
Suppose for $q<1$ and some $C>0$ we have
$ a_{n+1} < q a_n + C $, then $a_n$ is bounded by
$$a_n<\frac{C}{1-q}+a_0\;.$$
\end{lemma}

\subsection{Point forcing - Mean energy bounds}

Consider hyperviscous Burgers equation with point-forcing.
We would like to proceed exactly the way, we did in the previous section,
But we can not, as for point forcing, the solution of the  linear equation might fail to be in $L^\infty$.
This is why we add additional damping. Hyperviscous Burgers equation has been studied in several occasions.
See for example \cite{hyper3, hyper1, hyper2}.

Consider for some $\epsilon>0$ the operator $A_\epsilon=-\nu(-\Delta)^{1+\epsilon}$,
where $\Delta$ is the Laplacian subject to Dirichlet boundary conditions.
Then the hyperviscous Burgers equation is given by
\begin{equation}
  \partial_t u+ u\cdot \partial_x u  = A_\epsilon u +\alpha \delta_0 \dot\beta
\end{equation}
\begin{equation}
 u(\cdot,-L) = u(\cdot,L) =0, \qquad u(x,0)=u_0(x).
\end{equation}
Here, $\beta$ is a standard Brownian motion and $\delta_0$ the Delta-distribution.

Using the semigroup $e^{tA_\epsilon}$ the
solution for this system is (see \cite{DaPrato2} or
\cite{DaPratoB}):
\begin{equation} \label{pointsolu}
u(t) =  e^{tA_\epsilon}u_0 -\int_0^t e^{(t-\tau)A_\epsilon}(\lambda \Phi_\lambda(\tau)+\frac12\partial_x u(\tau,x)^2) d\tau +  \Phi_\lambda(t)
\end{equation}
where for some $\lambda\ge0$ fixed later
$$ \Phi_\lambda(t)= \alpha\int_0^t e^{(t-\tau)(A_\epsilon-\lambda)}\delta_0(x)  d\beta(\tau)
$$
solves
$$ \partial_t \Phi = A_\epsilon \Phi -\lambda \Phi+\alpha \delta_0 \dot\beta
$$
subject to
$$  \Phi(\cdot,-L) = \Phi(\cdot,L) =0, \qquad \Phi(x,0)=0\;.
$$
Using the standard orthonormal basis $\{e_k\}_{k\in\mathbb{N}}$ of eigenfunctions
of $A_\epsilon$ given  by $e_k(x)=\sqrt{1/L}\cdot\sin(-L+\frac{\pi k x}{2L})$ with corresponding eigenvalues
$\lambda_k=-(\pi k/2L)^{2+2\epsilon} $, we see
\begin{equation} \label{e:Foser}
\Phi_\lambda(t)= \alpha \sum_{k=1}^\infty \int_0^t e^{(t-\tau)(\lambda_k-\lambda)}d\beta(\tau) e_k(0)e_k\;.
\end{equation}
Note that the Fourier-coefficients of that series are not at all independent.
Thus we cannot rely on the better regularity results available for the stochastic
convolution of the previous chapter. Especially, for $\epsilon=0$ we cannot show that $\Phi_\lambda(t)\in L^\infty([-L,L])$.

Note that the series expansion of boundary and point forcing is very similar.
Thus we can regard boundary forcing at a point forcing at the boundary,
when the equation is subject to Neumann boundary conditions.

Our main result is now:
\begin{theorem}
For all $\epsilon>0$ and all initial conditions  $u_0$ independent of $\beta$   with $ E\|u_0\|^2 <\infty$
the solution of (\ref{pointsolu})
satisfies that the mean energy is uniformly bounded in time. I.e.,
$$ \sup_{t\ge 0}\EX \|u(t)-\overline{u}(t)\|^2 <\infty\;.
$$
\end{theorem}
We will proceed exactly as in the previous section.
Now  $v=u-\Phi_\lambda$ is a weak
solution of
\begin{equation}
\label{e:veqpoint}
  \partial_t v+ \frac12\partial_x (v + \Phi_\lambda)^2  = \nu \partial_x^2 v+\lambda \Phi_\lambda\;,
\end{equation}
again subject to Dirichlet boundary conditions and initial condition $v(0)=u_0.$

Now consider first the nonlinear term for some small $\delta>0$.
Using H\"older, Sobolev embedding of $H^{\frac12-\frac1p}$ into $L^p$ and
the bound
$$\|u\|_{H^{2\gamma}}\le C \|A_\epsilon^{\gamma/(1+\epsilon)}u\|$$
yields
\begin{eqnarray*}
 \int_{-L}^L v\Phi_\lambda v_x dx
&\le&\|v\|_{L^{2+\delta}} \|\Phi_\lambda\|_{L^{(4+2\delta)/\delta}} \|v_x\|\\
&\le& C\|A_\epsilon^{\frac14\frac{1}{1+\epsilon}\frac{\delta}{2+\delta}} v\|
\|\Phi_\lambda\|_{L^{(4+2\delta)/\delta}}
\|A_\epsilon^{\frac12\frac{1}{1+\epsilon}} v\|\;,
\end{eqnarray*}
Now we can easily find an $\delta>0$ sufficiently small such that there is a $p=p(\epsilon)\in(2,\infty)$
such that (using interpolation inequality)
$$ \int_{-L}^L v\Phi_\lambda v_x dx
\le C \|v\|_{L^2} \|\Phi_\lambda\|_{L^p}\|A_\epsilon^{\frac12} v\|\;.
$$
Now we can use the same proof as in the section before.
We only need that $\Phi_\lambda(t)\in L^\infty(0,L)$.
To be more precise, an easy calculation using the series expansion of (\ref{e:Foser})
shows that for any $\epsilon>0$
$$
\sup_{t\ge0} \mathbb{E} \|\Phi_\lambda(t)\|^2_{L^\infty}
\le C_\epsilon\sup_{t\ge0}\mathbb{E}\|\Phi_\lambda(t)\|^2_{H^{\frac{1+\epsilon}2}}
\to 0 \quad \mathrm{as\ } \lambda\to\infty\;.
$$
It is now straightforward to prove an analog to Lemma \ref{lem:Fernique}.
The remainder of the proof is analogous to the section before.

Let us remark that we could even simplify that proof a little bit,
by avoiding second order exponentials of $\Phi_\lambda$.
In that case we could work with $\lambda=0$

\subsection{Body forcing - Transient Behavior}

Let us focus on Burgers equation with body forcing. The results
for hyperviscous Burgers with point-forcing are completely
analogous. We will prove:
\begin{theorem}
\label{thm:tto0}
Let $u$ be a solution of (\ref{e:burgbodfo})
and consider for simplicity $u(0)=0$.
Denote by $$ E_u(t) =\EX\|u(t)-\overline{u}(t)\|^2$$ the mean energy of $u(t)$,
then there is some $\delta_0$  such that
$$
E_u(t) = E_{\Phi_0}(t) +\cO(t^{\frac12+\delta_0}) \quad\mathrm{\it for}\ t\to 0\;.
$$
To be more precise, for some $t_0>0$ sufficiently small there is a constant $C>0$
such that  $|E_u(t)-E_{\Phi_0}(t)| \le C t^{\frac12+\delta_0}$ for all $t\in[0,t_0]$.
\end{theorem}
As we know from  results like Theorem \ref{thm:mebo}
that $E_{\Phi_0}(t)$ behaves like $\sqrt{t}$ for small $t$,
we can conclude that the linear regime dominates for small $t$.

We could explicitly calculate $\delta_0$, but omit
this for simplicity of presentation.

For the proof of Theorem \ref{thm:tto0} use
$$|E_u(t)-E_{\Phi_0}(t)| \le C\EX\|v(t)\|^2\;,
$$
where we used Cauchy-Schwarz inequality and uniform bounds
on $\EX\|u(t)\|^2$ and $\EX\|\Phi_0(t)\|^2$.
Using (\ref{e:apri}) with $\lambda=0$ and $u(0)=0$ yields together with Lemma \ref{lem:Fernique}
$$\EX\|v(t)\|^2 \le C \int_0^t (\EX\|\Phi_0(t)\|^4_{H^{-1}})^{1/2} dt\;.
$$
It is now easy to show that $\EX\|\Phi_0(t)\|^4_{H^{-1}}$ behaves like $t^{2\delta_0}$ for some $\delta_0>0$,
which can be explicitly calculated using the methods of Theorem \ref{thm:mebo}.
Theorem \ref{thm:tto0} is now proved.

A simple corollary  using H\"olders inequality is:
\begin{cor}
Under the assumptions of Theorem \ref{thm:tto0},
we know for the mean correlation function
$$
C_u(t,r)=C_{\Phi_0}(t,r)+\cO(t^{\frac12+\delta_0}) \quad\mathrm{\it for}\ t\to 0 \ \mathrm{\it and\ all\ }r\;.
$$
\end{cor}
Notice that this result is only useful for small times and small $r$,
as seen from the qualitative behavior of $C_{\Phi_0}$, which is similar to
the results shown in section \ref{sec:corfu}.

\subsection{Trace class noise: Additive vs. multiplicative body noises}

Consider again a solution of the following Burgers equation:
\begin{equation}
  \partial_t u +u\cdot \partial_x u =\nu \partial_x^2 u + \sigma\dot{W}
\end{equation}
\begin{equation}
   u(\cdot,0) = 0,
  \qquad  u(\cdot,L) =0, \qquad u(x,0)=u_0(x),
\end{equation}
where $\{W(t)\}_{t\ge0}$ is a Brownian motion, with covariance $Q$, taking
values in the Hilbert space $L^2(0, L)$ with the usual scalar
product $\langle \cdot, \cdot \rangle$. We assume that the trace
$Tr (Q)$ is finite. So $\dot{W}$ is noise colored in space but
white in time.

Applying the It\^o's formula, we obtain
\begin{eqnarray}
   \frac12 d\|u\|^2 = \langle u, dW \rangle
   + [\langle u, u_{xx}-uu_x \rangle + \frac12 \sigma^2 Tr(Q)] dt.
\end{eqnarray}
as before $\langle u, uu_x \rangle=0$. Thus
\begin{eqnarray}
    \ddt \EX \|u\|^2
&=& -2 \|u_x\|^2  +  \sigma^2 Tr(Q).
\end{eqnarray}
By the Poincare inequality $\|u\|^2 \leq  c\|u_x\|^2$ for some
positive constant depending only on the length $L$, we have
\begin{eqnarray}
    \ddt \EX \|u\|^2  \leq  -\frac2{c} \|u\|^2   +  \sigma^2 Tr(Q).
\end{eqnarray}
Then using the Gronwall inequality, we finally get
\begin{eqnarray} \label{square-estimate}
      \EX \|u\|^2  \leq \EX \|u_0\|^2 e^{-\frac2{c}t}
        +  \frac12 c \sigma^2 Tr(Q) [1 - e^{-\frac2{c}t}].
\end{eqnarray}
Note that the first term in this estimate involves with initial
data, and the second term involves with the noise intensity
$\sigma$ as well as the trace of the noise covariance.

\bigskip

We now consider multiplicative body noise forcing.

\begin{equation}
  \partial_t u +u\cdot \partial_x u =\nu \partial_x^2 u + \sigma u
  \dot{w},
\end{equation}
with the same boundary condition and initial condition as above,
where $w_t$ is a scalar Brownian motion. So $\dot{w}$ is noise
homogeneous in space but white in time.

By the It\^o's formula, we obtain
\begin{eqnarray}
   \frac12 d\|u\|^2 = \langle u, \sigma u dw \rangle
   + [\langle u, \nu u_{xx}-uu_x \rangle + \frac12 \sigma^2 \|u\|^2] dt.
\end{eqnarray}
Thus
\begin{equation}
    \ddt \EX \|u\|^2
= -2 \nu \|u_x\|^2  +  \sigma^2  \|u\|^2  \nonumber \\
\leq  (\sigma^2-\frac{2\nu}{c}) \|u\|^2.
\end{equation}
Therefore,
\begin{eqnarray}
    \EX \|u\|^2  \leq  \EX \|u_0\|^2 e^{(\sigma^2-\frac{2\nu}{c}) t} .
\end{eqnarray}

Note here that the multiplicative noise affects the mean energy
growth or decay rate, while the additive noise affects the mean
energy upper bound.


\bigskip

{\bf Acknowledgments.}  Part of this work was done at the
Oberwolfach Mathematical Research Institute, Germany and the
Institute of Applied Mathematics, the Chinese Academy of Sciences,
Beijing, China. This work was partly supported by the NSF Grants
DMS-0209326
 \& DMS-0542450 and DFG Grant KON 613/2006.


\end{document}